\documentclass[aoas,preprint]{imsart}

\RequirePackage[OT1]{fontenc}
\RequirePackage{amsthm,amsmath}
\RequirePackage[numbers]{natbib}
\RequirePackage[colorlinks,citecolor=blue,urlcolor=blue]{hyperref}

\usepackage{wrapfig}
\usepackage{threeparttable} 
\usepackage{graphics}
\usepackage{graphicx}
\usepackage{amsthm}

\numberwithin{equation}{subsection}

\theoremstyle{plain}
\newtheorem{definition}{\textbf{Definition}}[section]
\theoremstyle{plain}

\theoremstyle{plain}
\newtheorem{remark}{\textbf{Remark}}[section]
\theoremstyle{plain}
\newtheorem{proposition}{\textbf{Proposition}}[section]
\theoremstyle{plain}
\newtheorem{lemma}{\textbf{Lemma}}[section]

\begin{document}

\begin{frontmatter}

\title{Object Oriented Data Analysis \\of Cell-Well Structured Data}
\runtitle{OODA of Cell-Well Structured Data}

\begin{aug}
\author{\fnms{Xiaosun} \snm{Lu}\thanksref{m1}\ead[label=e1]{xiaosun@live.unc.edu}},
\author{\fnms{J. S.} \snm{Marron}\thanksref{m1}\ead[label=e2]{marron@unc.edu}}
\and
\author{\fnms{Perry} \snm{Haaland}\thanksref{m1}\thanksref{m2}\ead[label=e3]{Perry\_D\_Haaland@bd.com}}

\runauthor{X. Lu et al}

\affiliation{University of North Carolina Chapel Hill\thanksmark{m1} and BD Technologies\thanksmark{m2}}

\address{University of North Carolina Chapel Hill\\
Department of Statistics and Operation Research\\
\printead{e1}\\
\phantom{E-mail:\ }\printead*{e2}}

\address{BD Technologies\\
21 Davis Drive, P.O. Box 12018\\
Research Triangle Park, NC 27709\\
\printead{e3}}
\end{aug}

\begin{abstract}
Object oriented data analysis (OODA) aims at statistically analyzing populations of complicated objects. This paper is motivated by a study of cell images in cell culture biology, which highlights a common critical issue: choice of data objects. Instead of conventionally treating either the individual cells or the wells (a container in which the cells are grown) as data objects, a new type of data object is proposed, that is the union of a well with its corresponding set of cells. This paper contains two parts.
The first part is the image data analysis, which suggests empirically that the cell-well unions can be a better choice of data objects than the cells or the wells alone. The second part discusses the benefit of choosing cell-well unions as data objects using an illustrative example and simulations. This research suggests that OODA is not simply a frame work for understanding the structure of the data analysis. It leads to useful interdisciplinary discussion that gives better results through more appropriate choice of data objects, especially for complex data analyses.

\end{abstract}

\begin{keyword}
\kwd{data objects}
\kwd{cell confluence}
\kwd{bright field image}
\end{keyword}

\end{frontmatter}

\section{Introduction}

The concept of Object Oriented Data Analysis (OODA) was introduced by Wang and Marron (2007) \cite{ooda}. The data objects are understood as the {\itshape atoms} of the statistical analysis. They could be numbers as taught in an elementary statistical course or vectors as in multivariate analysis. OODA, however, facilitates the analysis of populations of complex data objects. An interesting special case is functional data analysis, where the data objects are curves. See Ramsay and Silverman (2005) \cite{fda} for an overview of this type of analysis.  Dryden and Mardia (1998) \cite{shape} studied geometrical properties of objects, where the data objects are shapes. Wang and Marron (2007) \cite{ooda}, Aydin et al (2009) \cite{ooda2} and Shen et al (2013) \cite{shenfunctional} analyzed tree-structured data from medical images, where the data objects are trees.

Note that the concept of data objects generalizes the classical notion of that of experimental units.  An experimental unit is typically considered as one of a set of physical entities, each subjected to different experimental treatments. For instance, a {\itshape{well}} (i. e. a container for growing cells) with certain growth factors. On the other hand, OODA allows much more complex and abstract objects, such as images, shapes, trees, or even covariance matrices. 
The goal of OODA is to fully understand the data structure, choose appropriate data objects, and finally come up with an appropriate analysis oriented by this choice of data objects.
For example, in tree structured data analyses, combinatorial trees can be chosen as data objects to study tree structures.  In order to study the evolutionary relations among a group of organisms, phylogenetic trees are a good choice of data objects. See Holmes (1999 \cite{holmes1}, 2003a \cite{holmes2}, 2003b \cite{holmes3}) and Li et al (2000) \cite{ptree}. To exploit the power of functional data analysis to analyze data in tree space, the Dyck path representations are a good choice of data objects. See Shen (2012) \cite{dyck}. Note that OODA is about how to approach complex data analysis settings and is not limited to any particular data analysis methods. For example, nonparametric regression analysis of 3-d images as data objects was done in Davis (2008) \cite{davis} and of artery trees as data objects by Wang et al (2012) \cite{yuan}.

This paper discusses cell image analysis in cell culture biology from an object-oriented point of view. The motivation of this research is to develop a statistical approach to cell image analysis that better supports the automated development of stem cell growth media. A major hurdle in this process is the need for human expertise, based on studying cells under the microscope, to decide when to {\itshape{passage}} the cells to new media. We aim to use digital imaging technology coupled with statistical analysis to tackle this important problem (see Section \ref{SEC:background}). 
A new type of data object is proposed in Section \ref{SEC:ooda}. Comparison with other natural choices of data objects shows the benefit of this choice. Section \ref{SEC:objectLevelAnalysis} describes the final results of the image data analysis based on the choice of the proposed data object.
Section \ref{SEC:theory} further discusses the advantages of the proposed data object using an illustrative example and simulations, which can be easily generalized to any data set with a structure of groups and corresponding individuals. 

It is seen that OODA is not only a frame work for describing data objects, but also provides efficient terminology for making critical choices at the beginning of a complicated data analysis, especially in inter-disciplinary situations. In the example of this cell image analysis, biologists are comfortable with the notions of cell and well, but do not have simple terminology for the union. The discussion of ``what should we take as data objects?" allows quickly arriving at, and easy understanding by all parties involved of, the benefits of, the cell-well union as the best choice. 
Another excellent example of the benefits of OODA for facilitating inter-disciplinary discussion is in statistical acoustics research. See e.g. Aston et al (2012) \cite{sound}, where the raw data are digitally recorded sounds of human speech.  The data objects could be just the time series of sounds, but that might needlessly obscure key aspects of speech.  The data objects could also be any of various types of frequency analysis.  In the end, motivated by careful discussion of invariance principles, that interdisciplinary group finally chose a particular type of covariance matrices as data objects.

\subsection{Background of Cell Images}\label{SEC:background}

The maintenance and growth of cells under controlled conditions is called cell culture.  {\itshape{In vitro}} culture of cells taken directly from human tissues such as stem cells is, however, very difficult. Success depends on having the right conditions for growth, which include the type of container, the surface coating, oxygen levels, nutrients, and cell-signalling molecules. The liquid containing the nutrients and cell signalling molecules is generally called the growth medium. Two different growth media, having different components, may result in very different outcomes. There is great medical and commercial value in developing optimal growth media for stem cells, so the development of growth media is an important problem in the biotechnology industry. Furthermore, the use of automated methods to develop new media can greatly reduce development costs and increase the likelihood of success.

In order to produce enough cells for a medical procedure, cells are grown through several {\itshape{passages}} (or procedures). At each passage, cells are harvested and then reseeded into new vessels at a lower density, due to extensive cell-cell signalling as a function of density. Beyond a certain level of density, undesirable differences in morphology and phenotype arise (e.g. the cells are dying). So one of the most important problems in cell culture is deciding when to passage the cells. 
Cell density in a container is typically described in terms of {\itshape{confluence}}. The confluence of a cell culture is the percentage of the surface of the container that is covered by cells. For example, a 100\% confluent culture has cells in all surface area available for cell growth, whereas a 50\% confluent culture has used half of the available area. Usually it is desirable to passage a cell culture before it reaches 100\% confluence. In particular, stem cell cultures are often passaged at 80\% confluence.

\begin{figure}[htbp]
\centering
\includegraphics[scale=0.15]{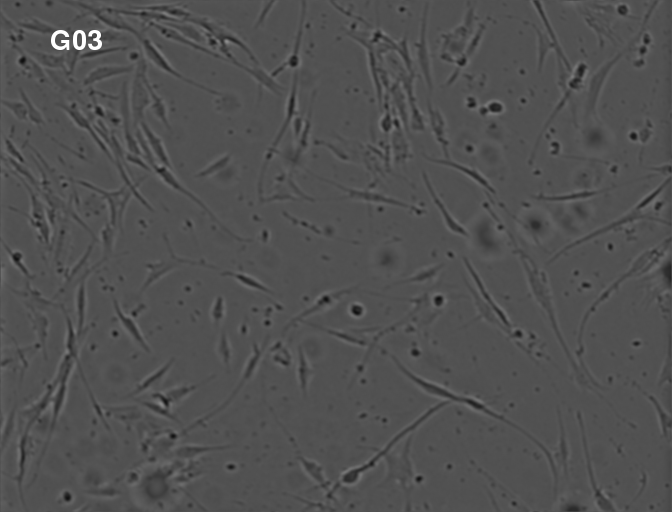} 
\includegraphics[scale=0.15]{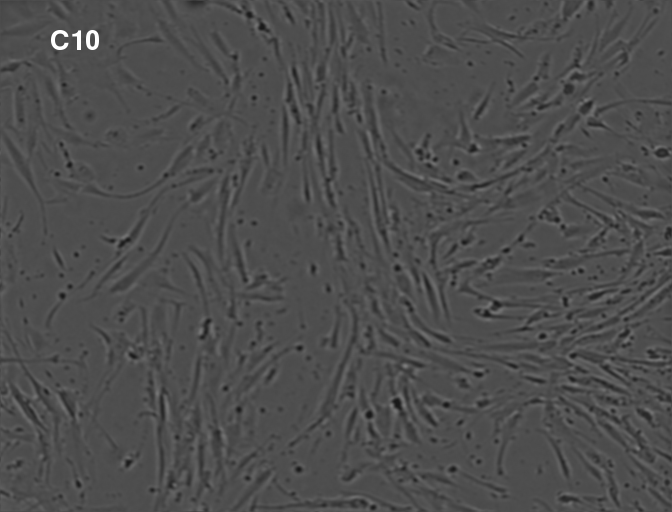} 
\includegraphics[scale=0.15]{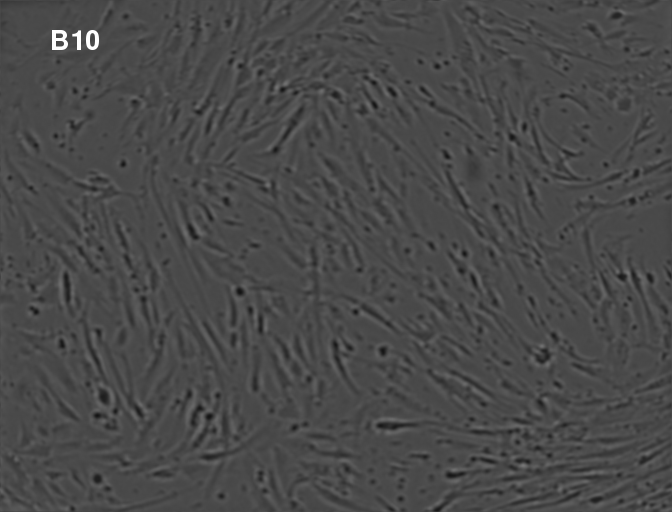}
\caption{Pre-processed and intensity-normalized bright field images
of three different wells from a 96-well plate of adherent stem cells, sorted from low
confluence level to high confluence level. The well names are on the upper
left corner. The cells correspond to the long thin objects. From
left to right, the cell number increases, the cell shape changes,
the gap between cells gets smaller and the cells begin to orient
with respect to each other.} \label{FIG:processedImages}
\end{figure}

Scientists often study images of the cells growing in the container to estimate the confluence and to decide whether or not to passage the cells. From a subjective viewpoint this is done by viewing the image and estimating the remaining space available for cell growth. This process is slow, manual, and highly variable, so being able to estimate the confluence directly from the image is an important capability of any automated cell growth platform. This estimation could be done, for example, by counting the number of cells, multiplying the number of cells by the average cell size and then comparing that area to the total surface area of the container. This approach is not generally desirable in an automated system because most methods to get this information kill the cells. A non-destructive way to get this information is through bright field imaging, where one shines a bright light down through the top of the container and records the image of the shadows from below. See Figure \ref{FIG:processedImages} for examples of bright field images. 

However, to determine cell confluence level based on the shadows in a bright field image is difficult. One can hardly tell the cell number in the image explicitly.
But some other visual factors in the
image can help biologists make their assessment of confluence level, such as the shape of the cells (more accurately, the cell shadows), the amount of empty space for the
cells to grow into and the {\itshape cell path} (the patterns in how
the cells orient with respect to each other). Changes of these visual factors as the confluence level increases can be seen in Figure \ref{FIG:processedImages}, where the three images are ordered from least confluent to most confluent. This manual assessment by biologists is usually subjective.
Thus, it is proposed to develop a statistical approach to numerically
summarize these visual features from an image and then make an objective statistical evaluation of cell confluence level.

\begin{figure}[htbp]
\centering
\includegraphics[scale=0.15]{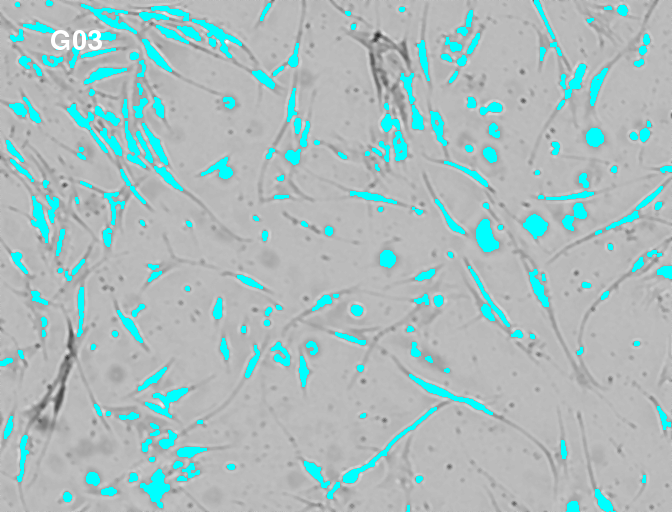} 
\includegraphics[scale=0.15]{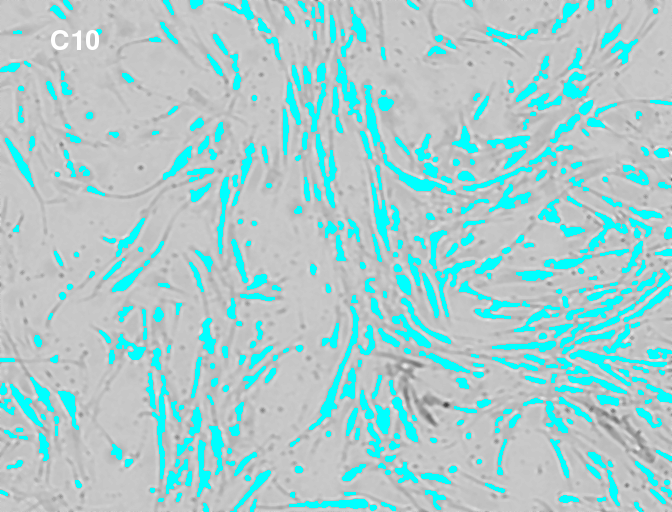} 
\includegraphics[scale=0.15]{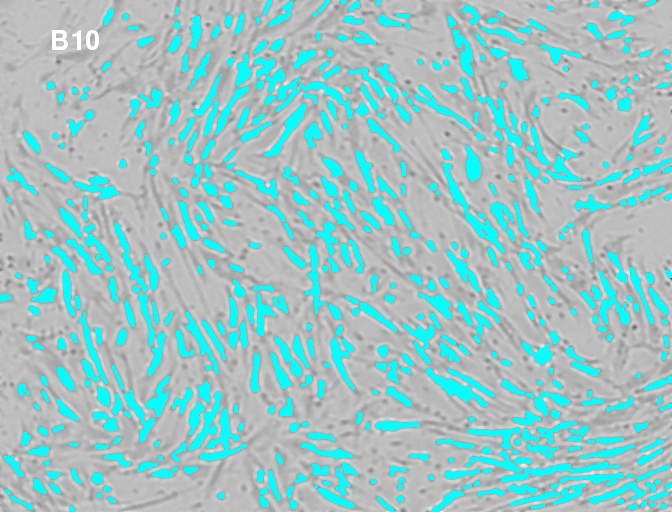}
\caption{Cell identification (using IPLab imaging software) of the wells shown in Figure \ref{FIG:processedImages}. The
well names are on the upper left corner. The cyan objects are the
identified cells.} \label{FIG:IPLabImages}
\end{figure}

A single 96-well plate of adherent stem cells from a screening
experiment by BD Technologies is selected as the training
sample. Each well is essentially a container in which cells are grown under a controlled condition. 
The culture conditions of the inner 60 wells represent a variety of culture
conditions that support different rates of cell growth, leading to
different confluence levels. The passaging decisions will be made on the well level, i.e. the cells in the same well will be passaged together. A bright field image is taken for each of the inner 60 wells (Figure \ref{FIG:processedImages}). The boundary of each cell is identified,
that is, the cell is segmented, using a custom script developed at BD
Technologies with IPLab for Pathway software (\url{http://www.digitalimagingsystems.co.uk/software/iplab.html}). Figure \ref{FIG:IPLabImages} shows the corresponding cell segmentation of the three bright field images in Figure \ref{FIG:processedImages}. Pixels that are
identified to be interior to cells are colored cyan. The identified objects do not
exactly cover the real cells, but this gives a useful approximation.

\begin{wrapfigure}{r}{0.41\textwidth}
\begin{center}
\includegraphics[scale=0.2]{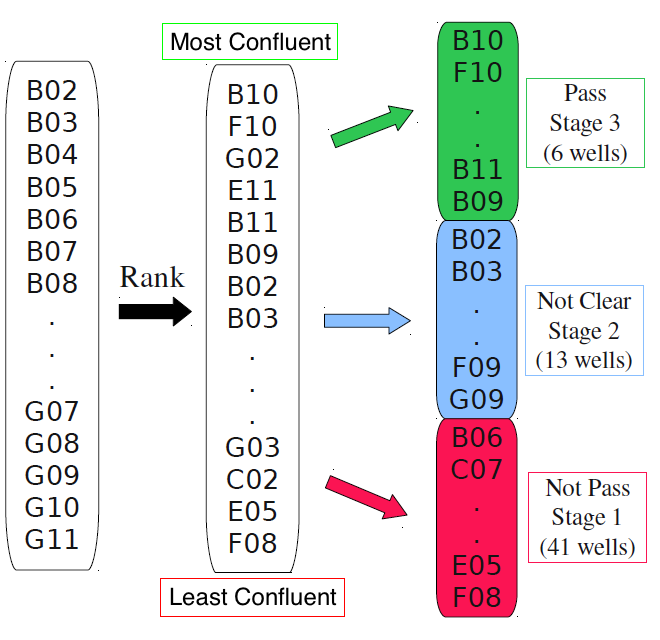}
\end{center}
\caption{Workflow of the manual assessment of confluence level. The
images were originally ordered by name (B02, B03, ...), and then
sorted in order of the estimated confluence by biologists. Finally the
passaging decisions were made based on the estimated confluence level: to passage (high
level), not clear (medium level) and not to passage (low level). }\label{FIG:rankProcess}
\end{wrapfigure}
Since the confluence level cannot be directly and unambiguously determined in a bright field image, in order to get a confluence evaluation of the 60 wells, an
experiment was designed where four biologists were asked to assess the confluence
level of the 60 images. 
Figure \ref{FIG:rankProcess} shows the
work flow of this experiment. The images were initially ordered by well
name, a random order of confluence level, as the condition of each well
was chosen under a randomized design. 
At first the biologists participated in the
experiment individually. Each of them sorted
the images in order based on their own estimated confluence level $\gamma$, and then specified two thresholds $\alpha_1$ and $\alpha_2$ ($\alpha_1>\alpha_2$) for making a
passaging decision for every image: to passage if $\gamma>\alpha_1$, not clear if $\alpha_2<\gamma<\alpha_1$, and
not to passage if $\gamma<\alpha_2$. 
However, the evaluation results
varied among biologists due to different subjective perceptions of confluence. 
After a careful
discussion, the biologists finally reached a consensus assessment, referred to later as {\itshape bio-assessment}, which will be considered as an unbiased evaluation of confluence level to judge the performance of the
statistical approach developed later in Section \ref{SEC:ooda}. 
This assessment resulted in each image receiving an integer indicating the {\itshape bio-rank} (the rank of confluence
level), and a categorical variable indicating the {\itshape bio-class}: low confluence level, medium confluence level or high confluence level. Each bio-class corresponds to a passaging decision: not to
passage, not clear, or to passage.

The goal of this research is to develop an objective and consistent approach for
assessing confluence level via statistical analysis of bright field images in order to better
support manual passaging decisions as well as provide the foundation
for an automated passaging system. The
conventional approach is {\itshape cell number assessment} (i.e. assessing cell confluence level
merely by counting the total number of the identified cells, ignoring other image
features), which does not match the bio-assessment very well. 
It is shown later in Section \ref{SEC:objectLevelAnalysis} that the alternative statistical approach proposed in this paper substantially improves the cell number assessment, in the sense of better predicting the bio-assessment.

\subsection{Feature Extraction}\label{SEC:features}

The bright field images are carefully pre-processed before extracting image features.
Some standard graphical techniques \cite{imageProcess}, such as flat field correction and convolution filter, are used to remove uneven background shading and granular noise. The intensity is normalized across images. 
Two types of confluence-related features are extracted from the images:
\begin{enumerate}
\item[(1)] Cell features, including properties of an individual cell and its relationship with its neighbors. These features can be categorized into four categories, intensity, shape \& size, local density, cell orientation (cell path), listed in Table \ref{TAB:summaryCellFea}.
\item[(2)] Entire-well features. Since cell confluence level is a function of the entire well instead of a
simple collection of cells, some additional well-level, or image-level, features are also considered in evaluating
confluence level.
These well-level features, such as the cell number and some summaries of the gap\footnote{The {\itshape gap} refers to the non-cyan area in IPLab segmented images (Figure
\ref{FIG:IPLabImages}). This gives an indication of how much more space the cells have to expand (generally, a smaller gap
indicates a higher confluence level). Also, as IPLab
identification of cells (the cyan objects) cannot exactly
cover the real cells, the gap contains part of the cell information.} 
in the image, are summarized in Table
\ref{TAB:summaryWellFea}.
\end{enumerate}
Due to irregular intensity distribution and irregular cell features respectively, two images are flagged as outliers.

\begin{table}[htbp]
\caption{Summary of cell features. }
\centering
\begin{tabular}{| c | p{2cm} | p{7.6cm} | c |}
\hline
& Categories    &  Details                                  & \# of Fea. \\
\hline\hline
1  & Intensity &  Average, Std., Average $log_{10}$, Minimum, Maximum, the 25\%, 50\%, 75\% quantiles of cell pixel intensity  
&8 \\
\hline 
2  &  Shape \& Size  & Perimeter, Area, Non-convexity, Length-Width Ratio, Radius Std.
& 5\\
\hline
3  &  Local Density  & Cell densities in 5 square moving windows with different sizes  & 5\\
\hline
4  &  Cell Orient.  & Cell angle, Angle difference with nearest neighbors, The 25\%, 50\%, 75\% quantiles of angle differences in 4 square moving windows with different sizes  
&14\\
\hline
\end{tabular}
\label{TAB:summaryCellFea}
\end{table}

\begin{table}[htbp]
\caption{Summary of additional entire-well features.  }
\centering
\begin{threeparttable} 
\begin{tabular}{| c | p{2cm} | p{7.6cm} | c |}
\hline
& Categories    &  Details                                  & \# of Fea. \\
\hline\hline
1  & Cell Number &     Number of identified cells in an image    &1 \\
\hline
2  &  Cell Gap  &  Summaries* of gap intensity 
 & 6\\
 \cline{3-4}
& & Summaries* 
of the size of circular gaps**& 6\\
\hline
\end{tabular}
\begin{tablenotes}
  \item[*] Standard deviation, min., max. and the 25\%, 50\%, 75\% quantiles are used as summaries.
  \item[**] These features are extracted by performing distance transformation \cite{distImage} on the IPLab segmented image. Statistical summaries of the intensity of the resulting distance image are used as a description of the size of the circular gaps among cells.
 \end{tablenotes}
\end{threeparttable}
\label{TAB:summaryWellFea}
\end{table}
\section{Object Oriented Data Analysis of Image Data}\label{SEC:ooda}

An important theme of OODA is that the very definition of data objects should be carefully considered before data analysis. In this cell image analysis, different choices of data objects are available and
lead to different results.
Since cell confluence level reflects the amount of available
space capacity of a well and the passaging decisions are made
at the well level, it is natural to treat wells as the data
objects. Meanwhile, as the cell features (Table \ref{TAB:summaryCellFea}) play an essential part in determining
confluence level, the individual cells should be considered as another
important aspect of the atoms of the analysis. 
Note that one could treat either the cells or wells alone
as data objects. Section \ref{SEC:comparison} shows a benefit from analyzing both the wells and the cells together, which motivates consideration of a new type of data object, that is the union of a well with its corresponding set of cells, or the {\itshape{cell-well unions}}. Section \ref{SEC:def} describes how the choice of data objects orients further analyses.

From an object-oriented point of view, the image data analysis is done in two steps:
(1) Separate analyses for various choices of data objects (Section \ref{SEC:comparison}), which show the advantage of treating the cell-well unions as data objects; 
(2) Analysis of cell-well unions as data objects (Section \ref{SEC:objectLevelAnalysis}), which provides the final results of our statistical assessment.
See Section \ref{SEC:theory} for further discussions of the choices of data objects.

\subsection{Data Objects and the Consequential Analyses}\label{SEC:def}
As discussed in Section \ref{SEC:features}, two different data sets are included in our analyses: 
\begin{itemize}
\item {\itshape{Cell data}} (containing cell features of each individual cell);
\item {\itshape{Well data}} (containing entire-well features of each well).
\end{itemize}
The cell sample size is always dramatically larger than the well sample size.
The first challenge in analyzing cell-well structured data is how to combine these two data sets. One natural solution is to define statistics to summarize the cell features across wells, and then combine the summarized cell data with the well data. Finally, the statistical passaging decision for each well will be made based on the combined data set. 

The following describes how the procedure of analysis will be oriented by the choice of data objects. Three different types of data objects and the corresponding data analyses are discussed.
\begin{enumerate}
\item[(1)]  {\itshape{Cells-alone analysis}}, i.e. analysis based on cells alone as data objects. In this analysis, only the cell data are used. The bio-assessment of a cell is defined the same as the bio-assessment of the well where the cell is cultured. The statistical passaging decision for a well is made based on the average of the predicted bio-classes of the individual cells in that well. Since all the additional entire-well features are ignored, one can expect that this analysis would not give a good classification of the passaging groups, i.e. the cells alone would not be an optimal choice of data objects.   
\item[(2)]  {\itshape{Wells-alone analysis}}, i.e. analysis based on wells alone as data objects. Both the cell data and the well data are used. However, since cells are not chosen as data objects, no cell data analysis is done here. The basic idea of this wells-alone analysis is to first summarize the cell features across wells directly by statistics, such as quantiles, and then combine the summarized cell data with the well data. Finally, the statistical passaging decisions are made by analyzing the combined data set.
\item[(3)]  {\itshape{Cell-well union analysis}}, i.e. analysis based on cell-well unions as data objects. This analysis uses both the cell data and the well data. First, the cell data analysis finds an appropriate way to summarize the cell data across wells. In particular, it finds a linear combination of the cell features that correlates well with the bio-rank, and then takes statistics, such as quantiles, of this linear combination and its orthogonal PC scores across wells as the summarized cell data. Finally, the statistical passaging decisions are made based on the combined data set of the summarized cell data and the well data. 
\end{enumerate}
\begin{figure}[htbp]
\centering
\includegraphics[scale=0.27]{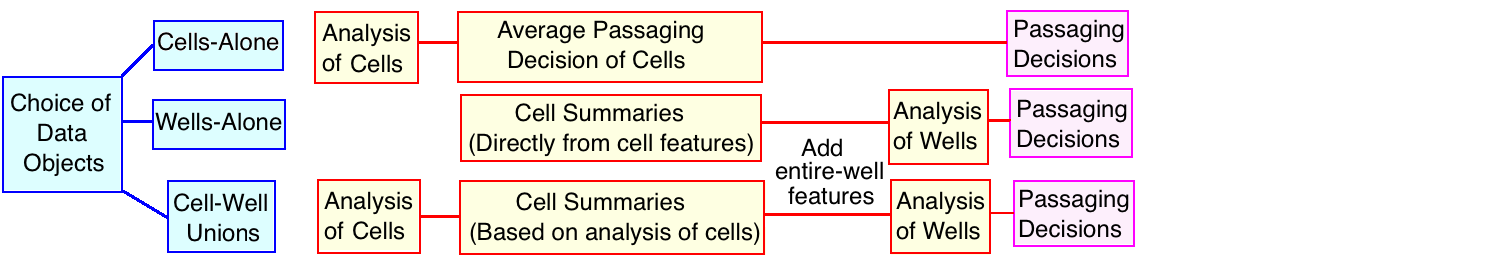}
\caption{Workflows of three different analyses, oriented by different choices of data objects respectively. Summaries refer to statistics, such as quantiles, of the cell-level features.} \label{FIG:dataObjects}
\end{figure}
The procedures of these three different analyses, oriented by the choice of data objects, are illustrated in Figure \ref{FIG:dataObjects}. Comparing the cell-well union analysis with the cells-alone analysis, it is seen that both of them begin with 
an analysis of the cell data. However, the cells-alone analysis makes passaging decisions simply based on the cell data analysis, while the cell-well union analysis includes an additional well-level analysis and makes passaging decisions based on the combined data set of the summarized cell data and the additional well data. It is also seen that the key difference between the cell-well union analysis and the wells-alone analysis is how the cell data are summarized across wells. The former incorporates an additional cell data analysis into the cell summarization. Further discussions of the choice of data objects in Section \ref{SEC:theory} concludes that the cell-well unions are a better choice of data objects than the other two choices.
\subsection{Comparison of Different Data Objects}\label{SEC:comparison}

This section aims at comparing the choices of three different data objects, the cells alone, the wells alone, and the cell-well unions, by performing the three corresponding analyses on the cell image data separately. 
For the purpose of comparison, we used the same statistical method, Distance Weighted Discrimination (DWD), to make the final passaging decisions in all these three analyses. 
Proposed by Marron et al (2007) \cite{dwd}, DWD is a powerful classification tool,  especially for high dimensional cases. It was used here to find the best linear separations between pairs of the three passaging groups and then to predict the group labels as the predicted passaging decisions. 
The consensus bio-classification, described in Section
\ref{SEC:background}, will be considered as a gold standard to judge the performance of these analyses.

\begin{enumerate}
\item[(1)]  Cells-alone analysis.  
Figure \ref{FIG:cellDataObj} (left) visualizes the cell data in two dimensions using Principal Component Analysis (PCA). The point color and the symbol are determined by the bio-assessment. The unclear pattern of either the colors or the symbols suggests that
the confluence information contained in the cell data is not obvious.
We intended to use DWD to classify the cell data directly. Each cell in a well would receive a label indicating its predicted passaging group, and the passaging decision for this well would be predicted by the average label of the cells within this well. However, due to the large sample size of the cells (over 20,000),  we encountered computational difficulties using the current DWD R package by H. Huang et al (2011)\cite{dwdR}. As an alternative approach, we randomly sampled the wells and randomly sampled a small set of cells from each well, and then used DWD to classify this smaller data set. This procedure was repeated 500 times, and the average classification error rate was 25.1\%. 
\item[(2)]  Wells-alone analysis.
Each cell feature was summarized into well-level features directly using 6 statistics: maximum, minimum, median, the 25\% and 75\% quantiles and standard deviation. The dimension of the summarized cell data is 6 times the original dimension. Then DWD was performed on the combined data set of the summarized cell data and the well data. The classification error rate was 8.6\%.
\item[(3)]  Cell-well union analysis. 
PCA was used to analyze the cell data, finding orthogonal directions that account for as much of the cell variability as possible. In Figure \ref{FIG:cellDataObj}, the left plot
shows a scatter plot of the first two PC scores, and the right shows only the averages across wells. It is seen that the vertical locations of the points reflect the order of the colors and the symbols, that is, PC1 reveals the bio-assessment. As a result of this PCA, each cell had totally 32 PC scores. The same 6-number summaries used in the wells-alone analysis were also used here to summarize the collections of scores across wells. The summarized PC scores, considered as the summarized cell data, were then combined with the well data. The DWD classification of this combined data set gives an error rate of 5.2\%. 
\begin{figure}[htbp]
\centering
\includegraphics[scale=0.44,page=1]{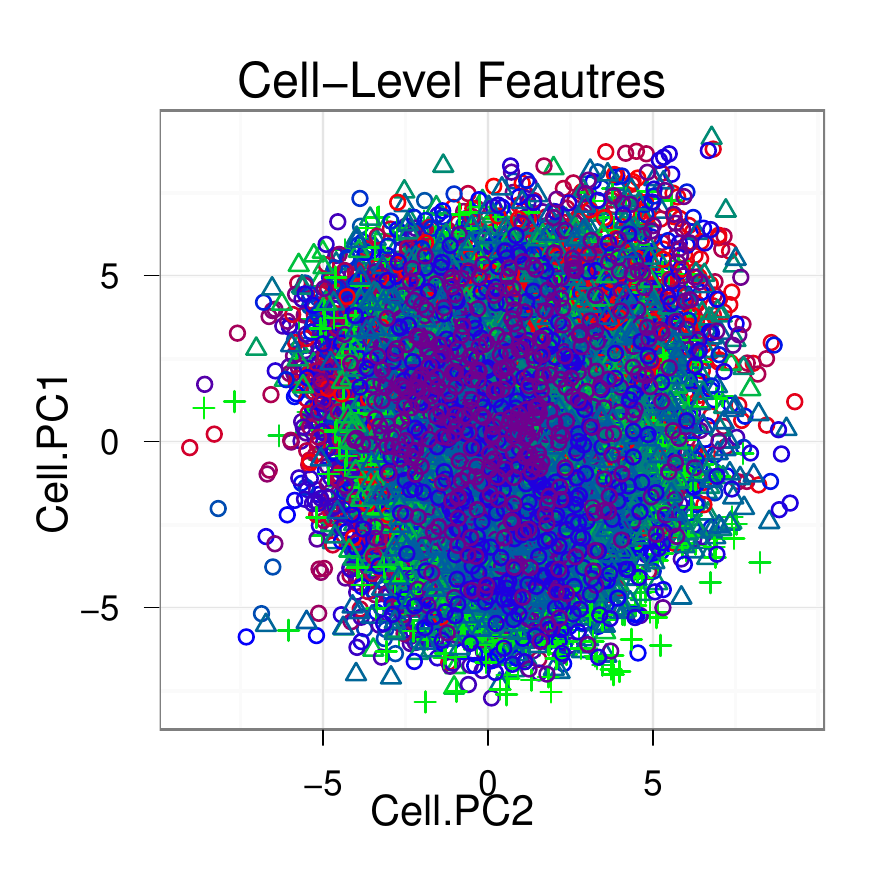}
\includegraphics[scale=0.44,page=2]{graph/cellDataObj.pdf}
\caption{PCA of the cell data. Left: Scatter plot of PC1 scores vs. PC2 scores. Right: Same plot as the left, only showing the averages across wells.
The points are colored from green (most confluent) to blue and then red (least confluent) according to the bio-rank. The symbols represent
the bio-classification: to passage (cross), not clear (triangle) and
not to passage (circle). PC1 conveys a lot of information about cell confluence. } \label{FIG:cellDataObj}
\end{figure}
\end{enumerate}

Compared with the cell-well union analysis, the cells-alone analysis has obvious disadvantages, as it ignores all the information in the well data and may also create computational challenges due to the large sample size of cells. In this image data study, the cell-well union analysis gives the lowest DWD classification error rate, and thus provides a set of statistical passaging decisions that is the most consistent with the bio-classification by biologists. A leave-one-well-out cross-validation shows the error rate of DWD classification using cell-well union analysis is 23\%, and wells-alone analysis 24\%.
This very slight advantage of the cell-well unions motivates a deeper look at this comparison.
Questions, such as how the cell-well unions gain an advantage over the wells alone, and whether the benefits from cell-well unions depend on statistical tools or the data structure, are discussed in Section \ref{SEC:theory}.

\subsection{Analysis of Cell-Well Data Objects}\label{SEC:objectLevelAnalysis}

Section \ref{SEC:comparison} suggests taking the cell-well unions as data objects and also describes the main procedure of the corresponding cell-well union analysis. This section provides the final results of the image data analysis as well as some supplementary details.

In the cell data PCA, the first four PCs totally explain nearly
70\% of the cell data variability.
Each of them reflects one of the four cell feature categories listed
in Table \ref{TAB:summaryCellFea}. Particularly,  PC1 is mainly about cell orientation, PC2 about cell intensity, PC3 about cell shape and size, and PC4 about local density.  
It is seen in Figure \ref{FIG:cellDataObj}  that the PC1 score correlates most to the bio-assessment.
Although most of the PCs do not correlate well with the bio-rank, we summarized all of the 32 PC scores across wells for further well-level analysis, in order to keep as much cell information as possible.
Experience suggests that dimension reduction may increase the error rate of predicting passaging groups and should be avoided.    

Finally, the percentage of false passaging decisions based on this cell-well union analysis, 5.2\%, is much lower than that from the DWD classification based on the cell number alone, 25.9\%
(A leave-one-well-out cross-validation shows that the former error rate is 23\%, while the latter is 30\%).
This result suggests that the statistical assessment based on image features can greatly improve the conventional cell number assessment, and thus can better support the automated passaging system.

\section{Illustrative Example and Simulations}\label{SEC:theory}
This section aims at exploring the potential generality of the superiority of cell-well union data objects using a toy example and simulations. 
As it is shown in Section \ref{SEC:ooda} that the cell-well union analysis has obvious advantages over the cells-alone analysis, this section will focus on comparing the cell-well union analysis with the wells-alone analysis. 

Figure \ref{FIG:dataObjects} shows that, in both the cell-well union analysis and the wells-alone analysis, one important step is to summarize the cell data across wells.
The essential difference between these two analyses is whether to summarize based on a cell data analysis or not. After cell summarization, there is no difference between the workflows of these two analyses. Hence, the following discussions will focus on the cell data analysis and cell summarization. 
Section \ref{SEC:toyExample} shows how the cell summarization can dramatically affect the result of the analysis using a two-dimensional toy example. Section \ref{SEC:summary} extends the toy example into more general cases, and concludes that 
the cell-well union summaries are generally better than the wells-alone summaries. Section \ref{SEC:simu} uses simulations to support the conclusion. 

In order to focus on the comparison of data objects, the following discussions are independent of any particular statistical tools that are used to analyze either the cells or the wells. In fact, the study of data objects provides suggestions of the choice of statistical tools as well as the choice of data objects. The basic idea is that, instead of comparing the final results from the cell-well union analysis and the wells-alone analysis, we compare the data patterns of the summarized cell data. Particularly, we assume that, if the summarized cell data in one analysis show a more clear pattern of the bio-assessment, referred to later as {\itshape{bio-pattern}}, then, no matter what statistical tools are used later to analyze this summarized cell data (or the combined data set of this summarized cell data and well data), it is easier to estimate the bio-assessment, and thus more probable to get a consistent estimation of the bio-assessment.


\subsection{Toy Example}\label{SEC:toyExample}
This section aims at illustrating how different ways of cell summarization lead to different well-level patterns of the summarized cell data, using a two-dimensional toy example. Let the cell data be $(x_{1},x_{2})$.
For convenient visualization, we only consider {\itshape{one dimensional cell summarization}} here, that is, each cell feature $x_{i}$ is summarized by a single statistic. After summarization, the dimension of the summarized well-level data is the same as the original cell data. 

\begin{figure}[htbp]
\centering
\includegraphics[scale=0.16,page = 1]{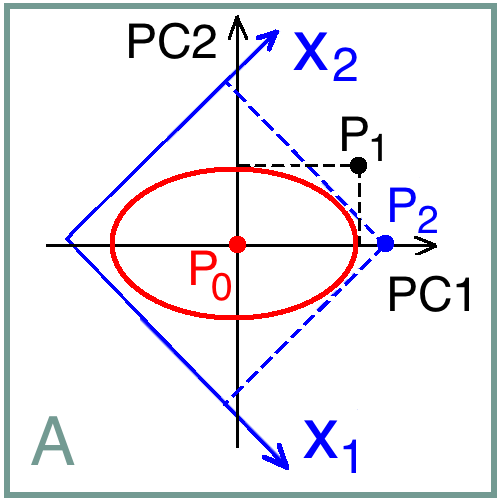}
\includegraphics[scale=0.15,page = 1]{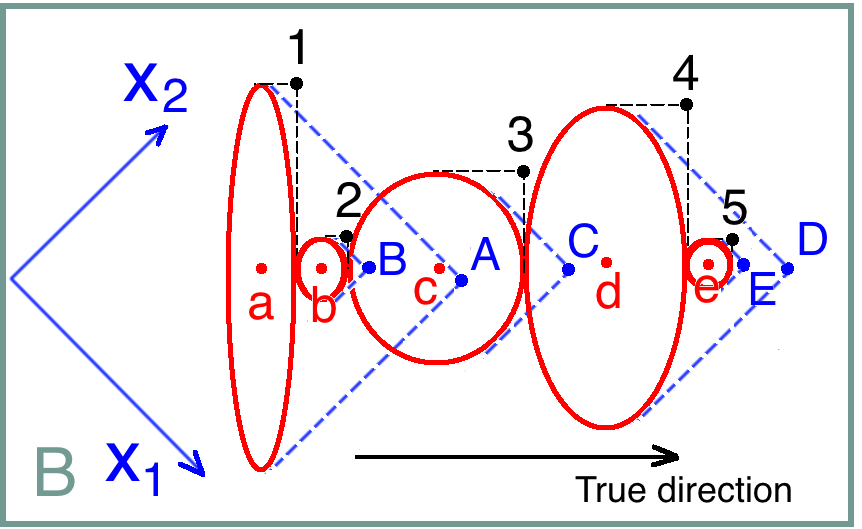}
\caption{Two-dimensional toy example based on maxima as summaries of the cells. Each red ellipse represents the distribution of the cell features in a well. Graph A: Two different summaries ($P_1$ and $P_2$) of the cell features of a single well, respectively corresponding to the cell-well union analysis (based on cell PCA, black) and the wells-alone analysis (based on the original cell features $x_1$ and $x_2$, blue). Graph B: How the different cell summarizations preserve or impair the underlying bio-pattern of cell data. The red points (a, b, c, d, e) are the population means of the wells, arranged along the true direction in order of bio-rank. The blue points (A, B, C, D, E) are the cell summaries based on cell features $(x_1,x_2)$, which result in poor estimates of the bio-rank. The black points (numbered 1, 2, 3, 4, 5) are based on the cell-level PCs, which give much better bio-rank estimates.  }\label{FIG:toyExample}
\end{figure}

Recall that, in the image data analysis in Section \ref{SEC:ooda}, the cell-well union analysis summarizes cell features based on their PC scores, while the wells-alone analysis is based on feature-wise summaries. This difference can be critical as highlighted by this toy example.  Figure \ref{FIG:toyExample}(A) illustrates the difference between the PC summarization and the direct summarization in a two-dimensional toy example based on only maxima as cell summaries. The red ellipse represents the cell feature distribution of a single well. The points $P_1$ and $P_2$ are the summarized well-level data from the two different cell summaries respectively. As long as $PC1$, $PC2$ are different from $x_{1}$, $x_{2}$, the corresponding summarized well-level data are different. The discussions here can be easily generalized to the cases of using other statistics, such as quantiles, standard deviation, etc. Note that taking the mean (the red point $P_0$) as a summary is a different case, where the summarized data from either the original cell features or their PCs are the same.

For the purpose of estimating the bio-assessment, we study the approaches to cell summarization for passing the bio-assessment information from cell data to further well-level analyses, that is, how the underlying bio-pattern in the cell data changes after cell summarization.     
Figure \ref{FIG:toyExample}(B) shows how the cell summarization can either impair or preserve the bio-pattern in cell data. 
Each of the five red ellipses represent the distribution of the cell features of a well. It is assumed that these wells have different bio-ranks, determined 
by their mean cell features (red points, labeled a, b, c, d, e, which are unknown in practice). 
The black arrow in the bottom right area shows the true direction of the bio-rank (practically unknown), which happens to be the same as the cell-level PC1. 
It is seen that the cell summaries based on ($x_1$, $x_2$), shown as the blue points (A, B, C, D, E), have a different order from the red ones (a, b, c, d, e), which lead to inconsistent estimates of the bio-rank. Thus the bio-pattern in the original cell data is impaired.
However, the summaries based on cell-level PCs, shown as the black points (numbered 1, 2, 3, 4, 5), give consistent estimates (i.e. the black numbers and the red lower case letters are in the same order). That is, the bio-pattern is well preserved after this cell PC summarization. 
Hence the cell-well union analysis is better than the wells-alone analysis. Note that the cell feature distributions of the wells vary a lot in this example. If those distributions are consistent, i.e. the shape and size of the red ellipses are all similar, one can imagine that the corresponding blue points (capital letters) and the black ones (numbers) will have the same order, i.e. the two sets of cell summaries will give the same bio-rank estimates.

In conclusion, how well the bio-pattern in the original cell data is preserved after cell summarization depends on both the summarizing method and the data structure. 
If the cell feature distributions vary across wells, then an additional cell-level analysis, such as PCA, can help construct a linear combination of the cell features that is close to the true direction and then better pass the bio-pattern in the cell data to further well-level studies by summarizing cells based on this linear combination. That is, the cell-well union analysis is better than the wells-alone analysis. 
On the other hand, if the cell feature distributions of the wells are consistent, then either cell summarization may perform equally well in capturing the bio-rank information from the cell data, that is, the cell-well union analysis and the wells-alone analysis may have equivalent performance.

\subsection{Cell Summarization}\label{SEC:summary}
This section extends the toy example into a more general case, and compares the cell-well union analysis with the wells-alone analysis by quantitatively studying how well the bio-pattern in cell data is preserved after cell summarization. Particularly, given that the cell feature means across wells reveal the underlying bio-rank, the variability of the bio-directional coefficient (defined later) of the summarized data provides a measurement of how well the bio-pattern is preserved. The main conclusions about the choice of data objects are in Remark \ref{REM:OODA}.

The following notations are used throughout this section. Consider $n$ wells and their $d_0$ dimensional cell features $X = (x_1, x_2,...,x_{d_0})$. 
Let $Y = (y_1, y_2,...,y_{d_0})$ be a set of orthogonal linear combinations of the cell features, and the statistics of Y across wells are taken as cell summaries. Note that $X$ is a special case of $Y$. In wells-alone analysis, $Y = X$. In cell-well union analysis, $Y$ can be, for example, the PC scores of $X$. Assume that the direction of the bio-rank exists and can be revealed by
the cell feature means $\mu = (\mu_1, \mu_2,...,\mu_{d_0})$.

First, for simplicity, consider one dimensional cell summarization, with each cell feature summarized by a single statistic. Let $\tilde{Y} = (\tilde{y}_1,\tilde{y}_2,...,\tilde{y}_{d_0})$ be the summarized cell features. Denote the true direction as $(\alpha_1, \alpha_2,...,\alpha_{d_0})$, where $\sum_{i=1}^{d_0}\alpha_{i}^2=1$. The projection of the population mean $\mu$ of a well on this direction, $\sum_{i=1}^{d_0}\alpha_{i}\mu_i$, reveals the bio-rank. Thus we assume a linear relationship between $\mu$ and the bio-rank. 

\begin{wrapfigure}{r}{0.28\textwidth}
\centering
\includegraphics[scale=0.2,page = 1]{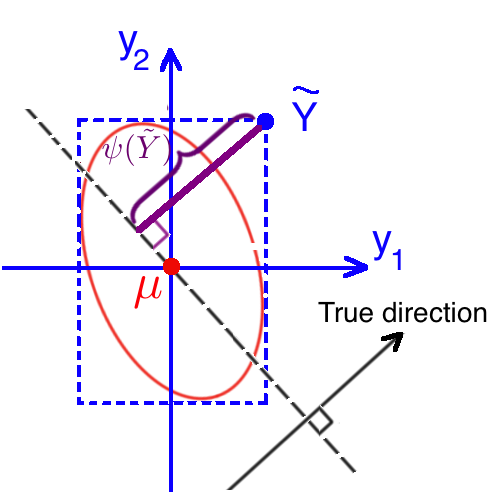}
\caption{The bio-directional coefficient $\psi(\tilde{Y})$ of the summary $\tilde{Y}$ (blue point) is shown as the purple line in a two dimensional case. The red point is the population mean $\mu$. }\label{FIG:dist}
\end{wrapfigure}

In order to measure how the bio-pattern is changed after cell summarization, for each well, we study the projection coefficient of the summarized point $\tilde{y}$ onto the true direction centered at the population mean $\mu$, referred to later as the {\itshape{bio-directional coefficient}}. This coefficient $\psi(\tilde{Y})$ is illustrated as the purple line in Figure \ref{FIG:dist}. 
Simple calculations show that 
\begin{equation}
    \psi(\tilde{Y}) = \sum_{i=1}^{d_0}(\tilde{y}_i-\mu_i)\alpha_{i}. 
    \label{EQU:dist}
\end{equation}
If the bio-directional coefficients of the wells are constant, then the bio-pattern is very well preserved after cell summarization. On the other hand, if these coefficients vary a lot, the bio-pattern is greatly impaired, or even lost, after cell summarization. Thus, the uncertainty of the preserved bio-pattern in the summarized cell data can be quantitatively expressed by the variability of these coefficients, defined as follows.

\begin{definition}\label{DEF: Uncertainty}
Consider a one dimensional cell summarization.
Let $\psi(\tilde{Y})$ be the bio-directional coefficient of the summarized data $\tilde{Y}$, defined in (\ref{EQU:dist}). Then the uncertainty of the bio-pattern in the summarized data, denoted as $\eta(\tilde{Y})$, is defined as the variance of $\psi(\tilde{Y})$, i.e.
$Var_{w}(\psi(\tilde{Y}))$, where the subscript ${w}$ highlights that it is a well-level variance.
\end{definition}

Then, assuming a Gaussian distribution of the cell features, we have the following lemma (see appendix for proofs). 
\begin{lemma}\label{LEM:1dUncertainty} 
Consider $d_0$ dimensional cell data $X$. The summarized data $\tilde{Y}= (\tilde{y}_1,\tilde{y}_2,...,\tilde{y}_{d_0})$ is derived by taking a single quantile of a collection of orthogonal linear combinations $Y = (y_1, y_2,...,y_{d_0})$ of the cell features.
Assume the cell data distributions of the wells are independent and Gaussian, with population mean $\mu = (\mu_1, \mu_2,...,\mu_{d_0})$. Assume the bio-rank direction $\alpha = (\alpha_1, \alpha_2,...,\alpha_{d_0})$ exists, where $\sum_{i=1}^{d_0}\alpha_{i}^2=1$, and is determined by $\mu$. Note that the $\alpha_i$'s depend on the coordinate system defined by $Y$, i.e. $\alpha = \alpha(Y)$.
Then the uncertainty of bio-pattern after cell summarization is
\begin{equation}
\eta(\tilde{Y}) = c_q^2 <\alpha^2(Y), Var_{w}Sd_{c}(Y)>,
\label{EQU:1dlemma}
\end{equation}
where $c_q$ is determined by the choice of the quantile, $<\cdot ,\cdot>$ denotes inner product, $Var_{w}Sd_{c}(Y)  = (Var_{w} Sd_{c}(y_1), ...\ , Var_{w} Sd_{c}(y_{d_0}))$
and the subscripts ${w}$ and ${c}$ indicate well-level and cell-level operations respectively.
\end{lemma}
Equation (\ref{EQU:1dlemma}) suggests that the uncertainty is bounded between \\
$c_q^2 \min_{i}\{Var_{w}Sd_{c}(y_i)\}$ and $c_q^2 \max_{i}\{Var_{w}Sd_{c}(y_i)\}$, regardless of the cell data dimension $d_0$.
Any cell summaries that lead to a smaller uncertainty will be considered better.
The uncertainty $\eta$ depends on the following three aspects.
\begin{enumerate}
\item[(1)]  The choice of the statistic, which is reflected by the term $c_q^2$ in the equation. Under the assumptions in the lemma, cell feature medians 
lead to a small $c_q$, and are the optimal choice. This is because the population mean of a well is assumed to determine its bio-rank. In practice, before cell summarization, it is always good to perform an exploratory analysis of the densities of $y_i$'s to choose a suitable quantile which nicely reflects the bio-pattern. Hence an additional cell-level analysis is always preferred. However, if $d_0$ is large, choosing proper quantiles for each $y_i$ is not feasible. 
\item[(2)] Well-level variability of the cell-level standard deviation, i.e. $Var_{w}Sd_c(Y)$. 
If the distributions of the cell-level data are consistent across wells, this term is 0. That is, whatever cell summaries are used, the uncertainty of the bio-pattern is 0. Thus both the wells-alone analysis and the cell-well union analysis give good estimates of the bio-rank. Under Gaussian assumptions, standardizing the cell-level data $Y$ across wells by their standard deviations can reduce this term. 
\item[(3)] The choice of the orthogonal linear combinations $y_i$, which is reflected by the term $\alpha^2(Y)$.
The inner product suggests that one should consider $\alpha(Y)$ and $Var_{w}Sd_{c}(Y)$ together. 
In the case of wells-alone analysis, $Y = X$, the inner product is $<\alpha^2(X), Var_{w}Sd_{c}(X)>$, which is determined by the structure of the original cell data .
In the cell-well union analysis, the cell-level analyses, such as PCA or Partial Least Square (PLS, taking the bio-rank as the response), can possibly construct a $Y$ that gives a smaller value of this inner product (See Section \ref{SEC:simu} for simulation results). Particularly, if $y_1$ captures the true direction, i.e. $y_1$ is the cell data projection on the true direction of the bio-rank,
then $\alpha_1=1$, $\alpha_k=0$ for $k\not=1$, thus the inner product is $Var_{w}Sd_{c}(y_1)$.
Assuming this true direction is reliable for estimating bio-rank in the sense that the cell data projections on it do not vary much across wells, this inner product can be very small. 
That is, the additional cell-level analysis can construct a better $Y$ to reduce the uncertainty of the bio-pattern after cell summarization.
Thus the cell-well union analysis can be better than the wells-alone analysis. Note that there should be no dimension reduction in $Y$ (e.g. all the PCs should be included when using PC summaries), because this may lose useful cell information for further well-level analysis.
\end{enumerate}

The above discussions can be easily generalized to multi-dimensional cell summarization. It is straightforward to extend Lemma \ref{LEM:1dUncertainty} to the following proposition.
The main lesson learned from the one dimensional cell summarization still holds.
\begin{proposition}\label{PRO:hdUncertainty} 
Consider a $d_s$ dimensional cell summarization, that is, the cell-level data $Y$ are summarized by $d_s$ quantiles across wells. The dimension of the summarized data $\tilde{Y}$ is $d_s d_0$. Let $\tilde{y}_{i j}$ be the $j$-th quantile of the original cell-level feature $y_i$, for $i=1,...,d_0$ and $j=1,...,d_s$. Suppose 
$\tilde{Y} = (\tilde{y}_{11},...,\tilde{y}_{1d_s},  ...,\tilde{y}_{d_0 1}...\tilde{y}_{d_0 d_s})$. Suppose the true direction in the summarized data space is of the form $(\alpha_{11}, ...,\alpha_{1d_s}, ..., \alpha_{d_0 1}, ...,\alpha_{d_0 d_s})$ where $\sum_{i=1}^{d_0}\sum_{j=1}^{d_s}\alpha_{ij}=1$. Under the same assumptions and notations of Lemma \ref{LEM:1dUncertainty}, the uncertainty of the bio-pattern after cell summarization is 
\begin{equation}
\eta(\tilde{Y}) = \sum_{j=1}^{d_s} c_{q(j)}^2 <\alpha_{(j)}^2(Y), Var_w Sd_c(Y)>,
\label{{EQU:hdlemma}}
\end{equation}
where $c_{q(j)}$ is determined by the choice of the $j$-th quantile, and
$\alpha_{(j)} = (\alpha_{1j}, ..., \alpha_{d_0 j})$. 
\end{proposition}
It is seen that the uncertainty $\eta(\tilde{Y})$ is bounded, regardless of $d_0$ and $d_s$, and its value also depends on the same three aspects as discussed previously in the case of one dimensional cell summarization.

As a conclusion, the cell-well union analyses are generally better than the wells-alone analyses, as stated in the following remark.
\begin{remark}\label{REM:OODA} 
Consider analyzing cell-well structured data. Assume the direction of the bio-rank exists.
In the process of summarizing the cell features, an additional cell-level analysis beforehand can help pass the underlying bio-pattern in the cell data to further well-level analysis more consistently. Thus the cell-well union analysis estimates the bio-rank better than the wells-alone analysis, or the cell-well unions are a better choice of data objects than the wells alone. The cases where the two types of data objects can be equally good are (1) The cell feature distributions across wells are consistent;
(2) One summarizing statistic is in the order of the bio-rank.
\end{remark}

Additionally, the previous discussion of Lemma \ref{EQU:1dlemma} suggests additional approaches to cell-level analyses. \begin{enumerate}
\item[(1)] Standardize the cell-level data of each well by their standard deviations, if the cell features are normally distributed;
\item[(2)] Choose statistics of the cell features that reflect the bio-rank, if feasible;
\item[(3)] Find a direction that nicely reveals the bio-rank (PCA and PLS are two recommended tools), and then summarize the cell features based on data projections in this direction and all its orthogonal directions.
\end{enumerate}

We investigated these approaches using the image data, including standardizing the cell data for each well and summarizing cell features using PLS instead of PCA. The standardization, however, did not improve the results much, because many cell features are not normally distributed and their quantiles can never be effectively standardized by the standard deviations. The PLS led to the same classification error rate as PCA, since the PLS direction (taking the bio-rank as the response) was very close to PC1. To further investigate these approaches, Section \ref{SEC:simu} uses simulations.

\subsection{Simulations}\label{SEC:simu}
This section validates the conclusions in Section \ref{SEC:summary} using simulations.

We simulated 50 wells, each with 50 to 300 cells, and each cell with 10 features. These cell features were normally distributed. The cell feature variances across wells were chosen as uniform (20, 500) random variables. The cell features of the wells had the same correlation structure (randomly generated). 
Cell data of different wells were simulated independently. The population means of the wells determined their bio-rank. These means were linearly located, and the difference between the means of two neighbor wells (with the bio-rank difference being 1) was 0.005. This bio-rank direction has equal entries, and thus has the same angle with each of the cell feature axes. The three passaging groups were defined by two thresholds on the population means. The data were then standardized. 

We did both wells-alone analyses and cell-well union analyses to estimate the passaging groups. Five quantiles, 1\%, 25\%, 50\%, 75\% and 99\%,  were used to summarize the cell-level data, and DWD was used to classify the groups. The wells-alone analyses were performed on two different cell-level data sets separately: the original cell data, and the cell data standardized within each well.
Two different cell summaries were considered in the cell-well union analyses: the PC summaries, and the PLS and its orthogonal PC summaries. 
Table \ref{TAB:simulations} shows the results of 500 simulations.
It is seen that additional cell analyses, such as PCA or PLS, reduce the classification error rate. Thus the cell-well union analyses are better than the wells-alone analyses. Also, lower uncertainty values lead to lower classification error rates, which is consistent with the discussions in Section \ref{SEC:summary}. 
Comparing the two wells-alone analyses suggests that standardizing cell data within each well reduces the classification error rate.

The simulations confirm that OODA can lead to a better choice of data objects, i.e. the cell-well unions, which leads to significantly better results than those from the wells-alone analyses.
\begin{table}[htbp]
\caption{Simulation results. }
\begin{threeparttable}
\begin{tabular}{| p{2.7cm}  |p{1.8cm} | p{1.8cm} | p{2.5cm} | p{2.39cm}|}
\hline
Data Objects & Wells-Alone   &  Wells-Alone    & Cell-Well Unions  & Cell-Well Unions \\
\hline
Cell Analyses & Not done & Std$^{[1]}$ & PCA \& Std$^{[1]}$ & PLS \& Std$^{[1]}$ \\
\hline\hline
 Uncertainty$^{[2]}$ &   $1.414 \pm 0.051$ &   $1.390 \pm 0.055$ & $0.471 \pm 0.088$ &$0.464 \pm 0.078$  \\
\hline
DWD Error Rate$^{[2]}$ &  $0.212 \pm 0.011$ & $0.132 \pm 0.009$  & $0.105 \pm 0.009$ & $0.104 \pm 0.009$\\
\hline
\end{tabular}
\begin{tablenotes}
  \item[1]   Standardize the cell data (or the PC/PLS scores) for each well by their standard deviations.
  \item[2]   The 95\% confidence intervals from 500 simulations are shown.
   \end{tablenotes}
\end{threeparttable}
\label{TAB:simulations}
\end{table}

\section{Conclusions}
This paper has proposed a new type of data objects, the cell-well unions, for the analysis of cell-well structured data, motivated by a study of cell images. We carefully discussed the choice of different data objects and compared their performances. It suggests that the cell-well unions are a better choice of data objects than either the wells alone or the cells alone. This paper clearly shows how the choice of data objects orients further analyses. 
In addition to just being a frame work for understanding the structure of the data analysis,
OODA, as effective terminology for inter-disciplinary communication, can guide critical choices of data objects, which can lead to better analyses of complex data.

\appendix

\section*{Acknowledgements}
We thank Warren Porter for helping with IPLab cell identification; Kirsten Crapnell, Tony Medlin, Divann Cofield and Richard Kelley for participating the experiment of assessing cell confluence level from images; and Tom Dyar and Dylan Wilson for the IT support.

\section*{Appendix} 
{\small
Proof of Lemma \ref{LEM:1dUncertainty}
\begin{eqnarray*}
\eta(\tilde{Y}) & = & Var_w(\psi(\tilde{Y}))\\
&=& Var_w\large\left(\sum_{i=1}^{d_0}(\tilde{y}_i-\mu_i)\alpha_{i}\large\right) 
           \hskip 0.3cm   (Equation\ \ref{EQU:dist})\\
           &=& \sum_{i=1}^{d_0}\alpha_{i}^2 Var_{w} (\tilde{y}_i-\mu_i) \hskip 0.3cm\\
&=& c_q^2\sum_{i=1}^{d_0}\alpha_{i}^2 Var_{w} Sd_{c}(y_i)\hskip 0.3cm (Under\ Gaussian\ assumption,\ \tilde{y}_i-\mu_i = c_q Sd_c(y_i))\\
\end{eqnarray*}
}

  \nocite{*}
\bibliographystyle{asa}
\bibliography{OODA-ArXiv_reference}

\end{document}